\documentclass[11pt]{article}

\usepackage{amsmath}
\usepackage{amssymb}
\usepackage{eucal}

\begin{document}

\baselineskip 16pt

\title{On  the intersection of the  $\cal F$-maximal subgroups  and
 the  generalized ${\cal F}$-hypercentre of a finite group }

\author{Wenbin Guo\thanks{Research of the first author is supported by
a NNSF grant of China (Grant \# 11071229)  and Wu Wen-Tsun Key Laboratory of Mathematics, USTC, Chinese Academy of Sciences.}\\
{\small Department of Mathematics, University of Science and
Technology of China,}\\ {\small Hefei 230026, P. R. China}\\
{\small E-mail:
wbguo@ustc.edu.cn}\\ \\
{ Alexander  N. Skiba \thanks{Research of the second author
 supported by Chinese Academy of Sciences Visiting Professorship for Senior
 International Scientists (grant No. 2010T2J12)}
}\\
{\small Department of Mathematics,  Francisk Skorina Gomel State University,}\\
{\small Gomel 246019, Belarus}\\
{\small E-mail: alexander.skiba49@gmail.com}}

\date{}
\maketitle

\begin{abstract}
Let $\cal F$ be a  class  of groups. A chief factor $H/K$ of a group
$G$ is called  \emph{${\cal F}$-central in $G$} provided
$(H/K)\rtimes (G/C_{G}(H/K)) \in {\cal F}$. We write $Z_{\pi{\cal F}}(G)$ to
denote  the product of all normal subgroups  of $G$ whose $G$-chief
factors  of order  divisible by at least one prime in $\pi$ are
$\cal F$-central. We call $Z_{\pi{\cal F}}(G)$ the $\pi{\cal
F}$-hypercentre of $G$.  A subgroup $U$ of a group $G$ is called
\emph{$\cal F$-maximal} in  $G$ provided that (a)  $U\in {\cal F}$,
and   (b) if $U\leq V\leq G$  and $V\in {\cal F}$, then $U=V$.
 In this paper we study the properties of the
intersection of all $\cal F$-maximal subgroups of a finite group. In
particular, we analyze the condition under which $Z_{\pi{\cal
F}}(G)$ coincides with the intersection of all $\cal F$-maximal
subgroups of $G$.
\end{abstract}


\let\thefootnote\thefootnoteorig

\section{Introduction}

Throughout this paper, all groups are finite and $G$ always denotes
a finite group. Moreover $p$ is always supposed to be a prime and
$\pi$  is a non-empty subset of the set $\Bbb{P}$ of all primes. We
use ${\cal G}_{\pi}$  (${\cal S}_{\pi}$) to denote the class of all
$\pi$-groups (of all soluble $\pi$-groups, respectively). In
particular, ${\cal G}_{p}$   denotes the class of all $p$-groups;
and we put that ${\cal G}_{\varnothing}={\cal S }_{\varnothing} =(1)$ is the class of all
identity groups.  We also use $\cal N$,  $\cal U$  and $\cal S$ to denote the
classes of all nilpotent groups,  of all supersoluble groups and  of all soluble groups,
respectively.

Let ${\cal F}$ be a class of groups. A group $G$ is said to be
 \emph{${\cal
F}$-critical} if $G$ is not
 in ${\cal F}$ but all proper subgroups of $G$
are in ${\cal F}$ \cite[p. 517]{DH}. If $1\in {\cal F}$, then we
write $G^{\cal F}$ to denote the intersection of all normal
subgroups $N$ of $G$ with  $G/N\in {\cal F}$. For any two classes $ {\cal
F}$ and ${\cal X}$ of groups, ${\cal X}{\cal F}$ is the class of
groups $G$ such that $G^{\cal F}\in {\cal X}$.

A \emph{formation} is a class ${\cal F}$ of groups with the
following properties: (i) Every homomorphic image of any group in
${\cal F}$ belongs to ${\cal F}$; (ii) If ${\cal F}\ne \varnothing$, then
 $G/G^{\cal F}\in {\cal F}$
for any group $G$. A formation ${\cal F}$ is said to be:
\emph{saturated} if $G\in {\cal F}$ whenever $G/\Phi (G) \in {\cal
F}$; \emph{hereditary} if $H\in {\cal F}$ whenever $H\leq G \in
{\cal F}$.

For any formation function $ f:\Bbb{P}\to \{ \text{formations of
groups}\}$, the symbol $LF(f)$ denotes the collection of all groups
$G$ such that either $G=1$ or $G\ne 1$ and $G/C_{G}(H/K)\in f(p)$
for every chief factor $H/K$ of $G$ and every $p\in \pi (H/K)$. If
${\cal F}=LF(f)$ for some formation function $f$, then $f$ is said
to be a local definition or a \emph{local satellite} (Shemetkov) of ${\cal
F}$. Every non-empty saturated formation ${\cal F}$ has a unique
local satellite $F$ with the following property: For any prime $p$,
both   $F(p)\subseteq  {\cal F}$  and $G\in F(p)$ whenever
$G/O_{p}(G)\in F(p)$ (see \cite[IV, Proposition 3.8]{DH}). Such a
satellite is called the \emph{canonical local satellite} of
${\cal F}$.

A chief factor $H/K$ of a group
$G$ is called  \emph{${\cal F}$-central in $G$} provided
$(H/K)\rtimes (G/C_{G}(H/K)) \in {\cal F}$. A normal subgroup $N$ of $G$ is said to be
\emph{$\pi{\cal F}$-hypercentral in $G$} if either $N=1$ or $N\ne 1$
and  every  chief factor of $G$  below $N$ of order  divisible by at
least one prime in $\pi$ is  ${\cal F}$-central in $G$. The symbol
$Z_{\pi{\cal F}}(G)$   denotes the \emph{$\pi{\cal F}$-hypercentre} of $G$,
that is, the product of all normal  $\pi{\cal F}$-hypercentral
subgroups of $G$. It is clear that for the  ${\cal F}$-hypercentre
$Z_{\cal F}(G)$ of $G$ ( see \cite[p.~389]{DH}) we have $Z_{\cal
F}(G)=Z_{\Bbb{P}{\cal F}}(G)$. On the other limited case, when $\pi
=\{p \}$,  $Z_{p{\cal F}}(G)$ is the the product of all normal
subgroups $N$  of $G$ such that every chief factor of $G$ below $N$
of order   divisible by  $p$ is  $\cal F$-central.

A subgroup $U$ of $G$ is called \emph{$\cal F$-maximal} in  $G$
provided that (a)  $U\in {\cal F}$,   and   (b) if $U\leq V\leq G$
and $V\in {\cal F}$, then $U=V$ \cite[p. 288]{DH}. We denote  the
intersection of all   $\cal F$-maximal subgroups of  $G$ by
$\text{Int}_{\cal  F}(G)$. In the paper \cite{BeidH}, Beidleman and
Heineken characterized the subgroup $\text{Int}_{\cal  F}(G)$ in the
case when $G$ is soluble and $\cal F$ is a hereditary saturated
formation.
 In this paper, as a development of some results in \cite{Skibaja5} and \cite{skjpaa}, we
find some new  properties and applications of the subgroup $\text{Int}_{\cal  F}(G)$.

Baer \cite{BaerI} proved  that $\text{Int}_{\cal  N}(G)$ coincides
with the hypercentre $Z_{\infty }(G)=Z_{\cal N}(G)$ of  $G$. Later, in
\cite{Sidorov}, Sidorov proved that if  ${\cal F}$   the class of
all soluble groups $G$ of nilpotent length  $l(G) \leq r \
 (r\in {\mathbb N})$, then for each  soluble group $G$, the equality $Z_{\cal F}(G)=\text{Int}_{\cal
F}(G)$  holds.  In the papers \cite{Skibaja5}   and \cite{skjpaa},
the analogous results were  obtained for the classes of all
$p$-decomposable groups and of all groups $G$ with $G'\leq F(G)$ in
the universe of all  groups. As one of our results in this paper, we
shall also prove that the intersection of all maximal
$p$-nilpotent subgroups of $G$ coincides with the subgroup
$Z_{p{\cal N}}(G)$. But in general, $Z_{\pi{\cal F}}(G) \ne
\text{Int}_{\cal  F}(G) $, even when  ${\cal F}$ is the class of all
supersoluble (all $p$-supersoluble, for any odd prime $p$)  groups
 and $G$ is soluble (see Theorem A and Remark 4.8 in Section 4).

{\bf Definition 1.1.}  Let $\cal X$ be a non-empty class of groups
and   ${\cal F}=LF(F)$ be a   saturated formation, where $F$ is the
canonical local satellite of ${\cal F}$.  We say that  $\cal F$
satisfies the \emph{$\pi$-boundary condition}  (\emph{the boundary
condition} if $\pi=\Bbb{P}$) \emph{in $\cal X$} if $G\in {\cal F}$
whenever $G\in {\cal X}$ and $G$ is an $F(p)$-critical group for
at least one  $p\in \pi $.

We say  that  $\cal F$
satisfies the  \emph{$\pi$-boundary condition}
 if $\cal F$ satisfies the  $\pi$-boundary condition
in the class of all groups.

If ${\cal F}$ is a non-empty  formation with $\pi ({\cal F})=
\varnothing$, then ${\cal F}=(1)$,  and therefore for any group $G$ we have $Z_{\cal
F}(G)= 1=\text{Int}_{\cal F}(G)$. In the other limited case, when
${\cal F}={\cal G}$ is the class of all groups, we have $Z_{\cal
F}(G)= G=\text{Int}_{\cal F}(G)$. Similarly, if ${\cal F}={\cal S}$, then  $Z_{\cal
F}(G)= G=\text{Int}_{\cal F}(G)$ for every soluble group $G$.

For the general case, we shall prove the following.

{\bf Theorem A.}   {\sl  Let ${\cal F}$  be a  hereditary saturated
formation with   $(1) \ne {\cal F}\ne {\cal G}$. Let $  \pi\subseteq
\pi (\cal F)$. Then  the equality $$Z_{\pi \cal
F}(G)=\text{Int}_{\cal F}(G)$$ holds  for each group $G$ if and only
if  ${\cal N}\subseteq {\cal F}={\cal G}_{\pi'}{\cal F}$ and   $\cal F$ satisfies the
$\pi$-boundary condition. }

Note that $N(p)={\cal G}_{p}$, where $N$ is the canonical local satellite of ${\cal N}$. Hence   
every $N(p)$-critical group has prime order. Therefore 
${\cal N}$ satisfies the
boundary condition and so the above-mentioned Baer's result is a first corollary of Theorem A.

{\bf Theorem B.}   {\sl  Let ${\cal F}$  be a  hereditary  saturated
formation of soluble groups with    $(1) \ne {\cal F} \ne {\cal S}$.  Let $ \pi\subseteq
\pi (\cal F)$.  Then  the equality
$$Z_{\pi \cal F}(G)= \text{Int}_{\cal F}(G)$$  holds for each soluble
group $G$ if and only if ${\cal N}\subseteq {\cal F}={\cal S}_{\pi'}  
  {\cal F}$ and
$\cal F$  satisfies the boundary condition in the class  of all
soluble   groups. }

If for some classes  ${\cal F}$ and ${\cal M}$  of groups  we have
${\cal F}\subseteq {\cal M}$, then every ${\cal F}$-maximal subgroup
of a group is contained in some its  ${\cal M}$-maximal subgroup.
Nevertheless, the following example shows that in general,
$\text{Int}_{\cal F}(G)\not\leq \text{Int}_{\cal M}(G)$.

{\bf Example 1.2.} Let ${\cal F}={\cal U}$  and ${\cal M}$ be the class of all
$p$-supersoluble groups, where $p > 2$. Let  $q$ be a prime dividing
$p-1$ and $G=P\rtimes (Q\rtimes C)$, where $C$ is a group of order
$p$, $Q$ is a  simple ${\mathbb F}_{q}G$-module which is faithful
for $C$ and  $P$  is a  simple ${\mathbb F}_{p}G$-module which is
faithful  for $Q\rtimes C$.  Then, clearly, $P= \text{Int}_{\cal
F}(G)$ and $\text{Int}_{\cal M}(G)=1$.

This example is a motivation for the following our result.

{\bf Theorem C.} {\sl Let  ${\cal F}\subseteq {\cal M}=LF(M)$
be   hereditary saturated formations with $\pi\subseteq  \pi ({\cal
F})$, where $M$ is  the canonical local satellites of ${\cal M}$.}

(a)  {\sl Suppose that ${\cal N}\subseteq {\cal M}={\cal G}_{\pi'}{\cal M}$ and
${\cal F}$ satisfies the  $\pi$-boundary condition in $\cal M$. Then
the inclusion $$\text{Int}_{\cal F}(G)\leq \text{Int}_{\cal M}(G)$$
holds for each group $G$}.

(b)  {\sl If  every  (soluble) $M(p)$-critical group 

belongs to $\cal F$ for every $p\in \pi$, then ${\cal N}\subseteq {\cal 
M}$ and 
$$\text{Int}_{\cal
F}(G)\leq Z_{\pi{\cal M}}(G)$$ for every  (soluble) group $G$.}

Recall that a subgroup $H$ of  a  group $G$ is said to be \emph{${\cal
F}$-subnormal} (in the sense of Kegel \cite{KegSubn}) or \emph{$K$-${\cal
F}$-subnormal} in $G$ (see p.236 in \cite{Bal-Ez}) if either $H=G$ or
there exists a chain of subgroups $$H=H_{0} <  H_{1} <     \ldots <
H_{t}= G $$ such that either $H_{i-1}$ is normal in $H_{i}$ or
$H_{i}/(H_{i-1})_{H_{i}}\in {\cal F}$ for all $i=1,  \ldots , t$.

For any group $G$, we write $\text{Int}^{*}_{\cal F}(G)$ to denote
the intersection of all  non-$K$-${\cal F}$-subnormal $\cal
F$-maximal   subgroups  of $G$. The following theorem shows that for
any hereditary saturated formation ${\cal F}$ with ${\cal N}\subseteq {\cal F}$,
 the intersection of
all  non-$K$-${\cal F}$-subnormal $\cal F$-maximal  subgroups of
a group $G$ coincides with $\text{Int}_{\cal
F}(G)$.

{\bf Theorem D.}  {\sl Let ${\cal F}$  be a   hereditary saturated
formation containing all nilpotent groups. Then  the equality
 $$\text{Int}^{*}_{\cal F}(G)=\text{Int}_{\cal
F}(G)$$  holds  for each group $G$.}

 We prove Theorems A, B, C and D  in Section 3. In Section 4 we give
some examples and discuss applications of these theorems.

All unexplained notation and terminology are standard. The reader is
referred to   \cite{DH},    \cite{Bal-Ez} and \cite{Guo} if
necessary.

\section{Preliminaries}

In view of Proposition 3.16 in \cite[IV]{DH},  we have

{\bf Lemma 2.1. } {\sl  Let ${\cal F}=LF(F)$  be a  hereditary  saturated
formation, where $F$ is the canonical local satellite of ${\cal F}$.
  Then for any prime $p$, the formation
$F(p)$ is hereditary.}

We shall need  in
our proofs a few facts about the  $\pi{\cal F}$-hypercentre.

{\bf Lemma 2.2. }  {\sl   Let ${\cal F}=LF(F)$  be a    saturated
formation, where $F$ is the canonical local satellite of ${\cal F}$. Let
     $  \pi\subseteq
\pi (\cal F)$ and $\sigma =\pi (\cal F)\diagdown \pi$.  
Let $N$ and $T$ be  normal subgroups of $G$, and $A\leq G$. }

(1) {\sl A chief factor $H/K$ of $G$ is ${\cal F}$-central if and only if $G/C_{G}(H/K)\in F(p)$
 for all primes $p$ dividing $|H/K|$.}

(2) {\sl Every  $G$-chief factor of $Z_{\pi{\cal F}}(G)$   of order
divisible by at least one  prime in $\pi$ is ${\cal F}$-central.}

(3) {\sl    $Z_{\pi{\cal F}}(G)N/N\leq  Z_{\pi{\cal F}}(G/N)$.}

(4) {\sl  $Z_{\pi{\cal F}}(A)N/N\leq  Z_{\pi{\cal F}}(AN/N)$.}

(5) {\sl If ${\cal F}$  is (normally) hereditary and $H$ is  a
(normal) subgroup of    $G$, then $Z_{\pi{\cal F}}(H)\cap E\leq
Z_{\pi{\cal F}}(H\cap E)$.}

(6) {\sl If  ${\cal G}_{\sigma}{\cal F}={\cal F}$ and $G/Z_{\pi \cal
F}(G)\in {\cal F}$, then $G\in {\cal F}$.}

(7) {\sl Suppose that  ${\cal F}$  is (normally) hereditary and let
$H$ be   a (normal) subgroup of    $G$. If  ${\cal G}_{\sigma}{\cal
F}={\cal F}$ and   $H\in {\cal F}$, then $Z_{\pi{\cal F}}(G)H\in
{\cal F}$.}

{\bf Proof.}  (1)  This assertion is well-known (see for example Theorem   
       17.14 in \cite{Shem-Sk} or  Theorem 3.1.6 in
\cite{Guo}). Assertions  (2) and (6) are evident.

(3) Let $H/K$ be a  chief factor of $G$ such that  $N\leq K < H\leq
NZ_{\pi{\cal F}}(G)$ and $|H/K|$ is divisible by at least one prime
in $\pi$. Then $H/K$ is $G$-isomorphic to the chief factor $H\cap
Z_{\pi{\cal F}}(G)/K\cap Z_{\pi{\cal F}}(G)$ of $G$. Therefore $H/K$
is  ${\cal F}$-central in $G$ by (1) and (2). Consequently, $Z_{\pi{\cal
F}}(G)N/N\leq Z_{\pi{\cal F}}(G/N)$.

(4) Let  $f:A/A\cap N \to AN/N$  be the canonical isomorphism from
$A/A\cap N $ onto $AN/N$. Then $f(Z_{\pi{\cal F}}(A/A\cap
N))=Z_{\pi{\cal F}}(AN/N)$ and
$$f(Z_{\pi{\cal F}}(A)(A\cap N)/(A\cap N))=Z_{\pi{\cal F}}(A)N/N.$$ By (3) we have $$Z_{\pi{\cal F}}(A)(A\cap N)/(A\cap N)\leq
Z_{\pi{\cal F}}(A/A\cap N).$$ Hence $Z_{\pi{\cal F}}(A)N/N\leq
Z_{\pi{\cal F}}(AN/N)$.

(5) First suppose that $\cal F$ is hereditary. Let $$1= Z_{0} < Z_{1}
< \ldots <  Z_{t} = Z_{\pi{\cal F}}(G)$$ be a chief
 series of $G$      below $Z_{\pi{\cal F}}(G)$ and  $C_{i}=
C_{G}(Z_{i}/Z_{i-1})$. Let $q$ be a prime divisor
 of $$|Z_{i}\cap H/Z_{i-1}\cap H|=|Z_{i-1}(Z_{i}\cap H)/Z_{i-1}|.$$ Suppose
 that $q$ divides  $|Z_{i}\cap H/Z_{i-1}\cap H|$.
 Then $q$ divide $|Z_{i}/Z_{i-1}|$, so $G/C_{i}\in F(q)$ by (1). Hence
 $H/H\cap C_{i}\simeq C_{i}H/C_{i}\in F(q)$. But  $H\cap C_{i}\leq
 C_{H}(Z_{i}\cap H/Z_{i-1}\cap H)$. Hence
 $H/C_{H}(Z_{i}\cap H/Z_{i-1}\cap H)\in F(q)$ for all primes $q$ dividing
$|Z_{i}\cap H/Z_{i-1}\cap H|$. Thus  $Z_{\pi{\cal F}}(G)\cap H\leq Z_{\pi{\cal F}}(H)$ by (1) and
 (2).
 But    then
 $$ Z_{\pi{\cal F}}(H) \cap E = Z_{\pi{\cal F}}(H) \cap
(H\cap E)\leq Z_{\pi{\cal F}}(H\cap E).$$ Similarly, one may prove the
second assertion of (5).

(7) Since $H\in { \cal F}$ we have  $$HZ_{\pi{\cal F}}(G)/Z_{\pi{\cal F}}(G)
\simeq H/H\cap Z_{\pi{\cal F}}(G) \in {\cal F}$$ and
 $$Z_{\pi{\cal F}}(G) \leq Z_{\pi{\cal F}}(Z_{\pi{\cal F}}(G)H)$$ by  (5).
 Hence $HZ_{\pi{\cal F}}(G) \in {\cal F}$  by (6).

The lemma is proved.

The   following lemma is evident (Note only that Statement (i) directly follows  from
 \cite[Theorem A.9.2(c), p. 30]{DH}.

{\bf Lemma 2.3.}  {\sl  Let ${\cal
F}$ be a  hereditary saturated formation.  Let $N\leq U\leq G$,
where $N$ is a normal subgroup of $G$. }

(i) {\sl If $G/N\in {\cal F}$ and $V$ is a minimal supplement of $N$ in $G$, then
 $V\in {\cal F}$.}

(ii) {\sl  If $U/N$ is an $\cal F$-maximal subgroup of $G/N$, then
$U=U_{0}N$ for some $\cal F$-maximal subgroup $U_{0}$ of $G$}.

(iii)  {\sl      If $V$ is an  $\cal F$-maximal subgroup of $U$,
then $V=H\cap U$ for some   $\cal F$-maximal   subgroup $H$ of $G$.}

The proofs of our theorems are based  on the following general facts
on the subgroup  $\text{Int}_{\cal F}(G)$.

{\bf Lemma 2.4.} {\sl Let ${\cal F}$  be a hereditary saturated
formation,     $  \pi\subseteq
\pi (\cal F)$ and $\sigma =\pi (\cal F)\diagdown \pi$.   Let $H$,
$E$ be  subgroups of $G$, $N$  a normal subgroup of $G$ and }
$I=\text{Int}_{\cal F}(G)$. {\sl Then: }

(a)   $\text{Int}_{\cal F}(H)N/N\leq \text{Int}_{\cal F}(HN/N)$.

(b)    $ \text{Int}_{\cal F}(H) \cap E\leq  \text{Int}_{\cal
F}(H\cap E)$.

(c)  {\sl  If $H/H\cap I\in {\cal F}$, then $H\in {\cal F}$. }

(d)  {\sl If $H\in {\cal F}$, then $IH\in {\cal F}$.}

(e) {\sl If $N\leq I$, then }
 $I/N= \text{Int}_{\cal F}(G/N)$.

(f)   $\text{Int}_{\cal F}(G/I)= 1$.

(g) {\sl  If ${\cal G}_{\sigma}{\cal F}={\cal F}$, then $Z_{\pi{\cal
F}}(G) \leq I$}.

{\bf Proof.} Assertions (a)-(f) are proved in \cite{Skibaja5}. Now
we prove (g).  Let $H$  be  a subgroup of $G$ such that  $H\in {
\cal F}$. Then  $HZ_{\pi{\cal F}}(G) /Z_{\pi{\cal F}}(G) \simeq
H/H\cap Z_{\pi{\cal F}}(G) \in {\cal F}$ and $Z_{\pi{\cal F}}(G)
\leq Z_{\pi{\cal F}}(HZ_{\pi{\cal F}}(G))$ by Lemma 2.2(5). It
follows from Lemma 2.2(5) that $HZ_{\pi \cal F}(G) \in {\cal F}$.
Thus $Z_{\pi{\cal F}}(G) \leq  I$.

The following lemma is a corollary of general  results on
$f$-hypercentral action (see \cite[Chapter IV, Section 6]{DH}). For
reader's convenience, we give a direct proof.

{\bf Lemma 2.5.} {\sl   Let ${\cal F}=LF(F)$  be a    saturated
formation, where $F$ is the canonical local satellite of ${\cal F}$.
Let $E$ be a normal $p$-subgroup of  $G$. If  $ E\leq Z_{\cal
F}(G)$, then $G/C_{G}(E)\in F(p)$.}

{\bf Proof.}  Let  $1=E_{0} < E_{1} <  \ldots  < E_{t} =E$   be a
chief series of $G$ below $E$. Let  $C_{i} =C_{G}(E_{i}/E_{i-1})$
and    $C= C_{1}\cap \ldots \cap C_{t}$. Then  $C_{G}(E)\leq C$ and
so $C/C_{G}(E)$ is a $p$-group by Corollary 3.3 in \cite[Chapter
5]{Gor}. On the other hand, by Lemma 2.2(1), $G/C_{i}\in F(p)$, so
$G/C\in F(p)$. Hence $G/C_{G}(E)\in F(p)={\cal G}_{p}F(p)$.  The
lemma is proved.

{\bf Lemma 2.6. }  {\sl Let  ${\cal F}=LF(F)$ and ${\cal M}$ be saturated formations with
 $p\in \pi ({\cal F})$ and ${\cal F}\subseteq {\cal M}$,
 where $F$ is the canonical local satellite of ${\cal F}$. Suppose that $G$ is a group
 of minimal order
in the set of all $F(p)$-critical  groups $G\in {\cal M}$ with
$G\not \in {\cal F}$. Then  $O_{p}(G)=1=\Phi (G)$ and $G^{\cal F}$
is the unique minimal normal subgroup of $G$. }

{\bf Proof.} Let $N$ be a minimal normal subgroup of $G$. First we show that  $G/N\in
{\cal F}\cap {\cal M}$. Indeed, since $G\in {\cal M}$ and ${\cal M}$ is a formation,
 $G/N\in {\cal M}$. Suppose that    $G/N\not \in     {\cal F}$. Then $G/N
\not \in F(p)  $ since  $F(p)\subseteq {\cal F}$.
 On the other hand, for any maximal subgroup $M/N$ of $G/N$
we have $M/N\in F(p)$  since $F(p)$ is a formation  and $G$ is an $F(p)$-critical
 group.   Thus  $G/N$  is an
  $F(p)$-critical  group in  with $G/N\not \in   {\cal M}\setminus   {\cal F}$, which
 contradicts the minimality of $G$.
 Hence  $G/N \in     {\cal F}$.
Since ${\cal F}$ is a saturated formation,  $N=G^{\cal F}$ is a unique
 minimal normal   subgroup of $G$ and  $\Phi (G)=1$. Suppose that  $N\leq  O_{p}(G)$  and
let $M$ be a maximal subgroup of $G$ such that $G=NM$. Then $G/N\simeq
M/N\cap M \in F(p)={\cal G}_{p}F(p)$, so $G\leq F(p)\subseteq {\cal F}$.
This contradiction shows that $O_{p}(G)=1$. The lemma is proved.

{\bf Lemma 2.7.} {\sl Let $\cal F$ be a formation, $H$ and $E$ be
subgroups of a group $G$, where $H$ is $K$-${\cal F}$-subnormal in
$G$. Then: }

(i) {\sl  $H\cap E$ is $K$-${\cal F}$-subnormal in   $E$ }(see Theorem 6.1.7 in \cite{Bal-Ez}).

(ii) {\sl If $E$ is normal in $G$, then $HE/E$ is   $K$-${\cal F}$-subnormal in
  $G/E$} (see Theorem 6.1.6 in \cite{Bal-Ez}).

{\bf Lemma 2.8. } {\sl  Let ${\cal F}$  be a  hereditary saturated
formation.  Let $N\leq U\leq G$, where $N$ is a normal subgroup of
$G$. }

(i) {\sl  If $U/N$ is  a non-$K$-${\cal F}$-subnormal  $\cal
F$-maximal subgroup of $G/N$, then  $U=U_{0}N$ for some
non-$K$-${\cal F}$-subnormal $\cal F$-maximal subgroup $U_{0}$ of
$G$}.

(ii)  {\sl      If $V$ is a non-$K$-${\cal F}$-subnormal  $\cal
F$-maximal subgroup of $U$, then $V=H\cap U$ for some non-$K$-${\cal
F}$-subnormal $\cal F$-maximal   subgroup $H$ of $G$.}

{\bf Proof.}  (i) By Lemma 2.3 (ii), there is an  $\cal F$-maximal
subgroup $U_{0}$ of $G$ such that $U=U_{0}N$. Since  $U/N$ is
non-$K$-${\cal F}$-subnormal in $G/N$, $U_{0}$ is  not  $K$-${\cal
F}$-subnormal in $G$ by Lemma 2.7(ii).

(ii) By Lemma 2.3(iii), for  some $\cal F$-maximal   subgroup $H$ of
$G$ we have  $V=H\cap U$. Since $V$ is a non-$K$-${\cal
F}$-subnormal in $U$, $H$ is  not  $K$-${\cal F}$-subnormal in $G$
by Lemma 2.7(i).

{\bf Lemma 2.9.}   {\sl Let ${\cal F}=LF(F)$  be a  non-empty
saturated formation, where $F$ is the canonical local satellite of
${\cal F}$.}

(1) {\sl  If  ${\cal F}={\cal G}_{p}{\cal F}$ for some prime $p$,
then  $F(p)={\cal F}$. }

(2) {\sl If  ${\cal F}={\cal N}{\cal H}$ for some non-empty
formation ${\cal H}$, then $F(p)= {\cal G}_{p}{\cal H}$ for all
primes $p$. }

{\bf Proof.} (1)  Since $ F(p)\subseteq {\cal F}$,  we need only
prove that $ {\cal F}\subseteq F(p)$. Suppose that this is false and
let $A$ be a group of minimal order in ${\cal F}\diagdown F(p)$.
Then $A^{F(p)}$ is a unique minimal normal subgroup of $A$ since
$F(p)$ is a formation.  Moreover,  $O_{p}(A)=1$ since $F(p)= {\cal
G}_{p}F(p)$. Let  $G=C_{p}\wr A= K\rtimes$ where $K$ is the base group of 
the regular wreath $G$. Then $K=O_{p', p}(G)$ and $G\in {\cal F}={\cal 
G}_{p}{\cal F}$. Hence         $A
\simeq G/P=G/O_{p', p}(G) \in F(p)$, a contradiction. Thus
$F(p)={\cal F}$.

(2) The inclusion $F(p)\subseteq {\cal G}_{p}{\cal H}$ is evident.  The  
inverse
 inclusion can be proved similarly as the inclusion  $ {\cal F}\subseteq F(p)$ in the proof of 
(1).

We will also use in our proofs the following well-known  elementary fact
 (see for example, \cite[Lemma 18.8]{Shem-Sk} or \cite[Lemma 3.5.13]{Guo}).  

{\bf Lemma 2.10.}  {\sl  
If  $O_{p}(G)=1$    and $G$ has a unique minimal normal subgroup, then there exists  a  simple 
   ${\mathbb F}_{p}G$-module
 which is     faithful  for $G$.  }

\section{Proofs of the Theorems}

{\bf Proof of Theorem C.}  (a) Suppose that this assertion   is
false  and let $G$ be a counterexample of minimal order.   Let $I=
\text{Int}_{\cal F}(G)$ and   $I_{1}= \text{Int}_{\cal M}(G)$.
 Then  $1 < I < G$ and  $I_{1} \ne G$.
Let  $L$ be  a minimal normal subgroup of $G$ contained   in $ I$ and $C=C_{G}(L)$.
Then $\pi (L)\subseteq \pi ( {\cal F})$.

(1)    $IN/N\leq \text{Int}_{\cal F}(G/N) \leq  \text{Int}_{\cal
M}(G/N)$  {\sl for any non-identity  normal subgroup $N$ of $G$.}

Indeed, by Lemma 2.4(a), we have $IN/N\leq  \text{Int}_{\cal
F}(G/N)$. On the other hand, by the choice of $G$,   $
\text{Int}_{\cal F}(G/N) \leq  \text{Int}_{\cal M}(G/N)$.

(2) {\sl $L\nleq I_{1}$; in particular, the order of $L$ is
divisible by some prime $p\in \pi$}.

Suppose that $L\leq I_{1}$. Then    $I_{1}/L= \text{Int}_{\cal
M}(G/L)$ by  Lemma 2.4(e).  But by (1), $IL/L=I/L \leq \text{Int}_{\cal
F}(G/L) \leq  \text{Int}_{\cal M}(G/L)$. Hence  $I/L\leq I_{1}/L$
an so $I\leq I_{1}$, a contradiction.  Thus $L\not\leq I_1$. This
means that there exists and ${\cal M}$-maximal subgroup $M$ of $G$
such that $L\nleq M$. Suppose that $L$ is a $\pi'$-group. Then $LM\in {\cal
G}_{\pi'}{\cal M}={\cal M}$, which contradicts the maximality of $M$.
Hence  the order of $L$ is
divisible by some prime $p\in \pi$.

(3) {\sl If $L\leq M < G$, then }$L\leq \text{Int}_{\cal M}(M)$.

By Lemma 2.4(b), $L\leq I\cap M\leq \text{Int}_{\cal F}(M)$. But
since $|M| < |G|$,  $ \text{Int}_{\cal F}(M)\leq  \text{Int}_{\cal
M}(M)$  by the choice of $G$.   Hence  $L\leq \text{Int}_{\cal
M}(M)$.

(4) {\sl $G=L U$ for any ${\cal M}$-maximal subgroup
$U$ of $G$  not containing $L$. In particular, $G/L\in {\cal M}$.}

Indeed, suppose that  $LU\ne G$.  Then by (3), $L\leq
\text{Int}_{\cal M}(LU)$, which implies that  $LU\in {\cal M}$ by
Lemma 2.4(c). This contradicts the ${\cal M}$-maximality of $U$.
Hence we have (4).

(5) {\sl $C_{G}(L)\cap U=U_{G}=1$ for any ${\cal M}$-maximal subgroup
$U$ of $G$  not containing $L$}.

Since $C_{G}(L)$ is normal in $G$ and $G=LU$ by (4),
$U_{G}=C_{G}(L)\cap U$. Assume that $U_{G}\ne 1$. Let $U/U_{G}\leq
W/U_{G}$, where $W/U_{G} $ is an ${\cal M}$-maximal subgroup of
$G/U_{G}$. Then by (1), $LU_{G}/U_{G}\leq W/U_{G}$. Hence $G=LU\leq
W$ by (4), which  means that $G/U_{G}=W/U_{G}\in {\cal M}$. But by
(4), $G/L\in {\cal M}$.  Therefore $G\simeq G/L\cap U_{G}\in {\cal
M}$, and consequently $I=G$, a contradiction. Hence (5) holds.

{\sl The final contradiction for (a).}

Since $L\nleq I_{1}$ by (2), there is an ${\cal M}$-maximal subgroup
$M$ of $G$ such that $L\nleq M$.  But then $G=L M$ by (4).  Since
$L\leq I$  and $G\not \in {\cal F}$, $M\not \in {\cal F}$  by Lemma
2.4(d). 
Let $H$ be an $\cal F$-critical subgroup of $M$, $V$ a maximal
subgroup of $H$. We show that $V\in F(p)$. By Lemma 2.4(d), $D=LV\in
{\cal F}$. Hence $D/O_{p',p}(D)\in F(p)$. First assume that $L$ is a
non-abelian group. Then, since $p$ divides $|L|$, $O_{p',p}(D)\cap
L=1$. Hence $O_{p',p}(D)\leq C_{G}(L)$ and $O_{p',p}(D)\cap V=1$ by
(5). Since ${\cal F}$ is hereditary, $F(p)$ is hereditary by Lemma
2.1. Therefore $O_{p',p}(D)V/O_{p',p}(D)\simeq V\in F(p)$. Now assume
that $L$ is an abelian $p$-group. Then $L\leq O_{p',p}(D)$ and
$O_{p',p}(D)=L(O_{p',p}(D)\cap V)$. Hence $ O_{p'}(D)\leq M\cap
C_{G}(L)=1$. It follows that $O_{p',p}(D)=O_{p}(D)$. Therefore
$D/O_{p}(D)\in F(p)={\cal G}_{p}F(p)$, which implies that $D\in
F(p)$ and so $V\in F(p)$. Therefore $H$ is an $F(p)$-critical
group. Since $\cal M$ is hereditary and $M\in {\cal M}$, $H\in {\cal
M}$.  But then $H\in {\cal F}$ since ${\cal F}$ satisfies the
$\pi$-boundary condition in $\cal M$ by hypothesis. This
contradiction completes the proof of (a).

(b) Suppose that    every   $M(p)$-critical group $G$
belongs to $\cal F$ for every $p\in \pi$.  First  we show that ${\cal N}\subseteq {\cal 
M}$. Assume that this is false and let $C_{q}$ be a group of prime order $q$ with  
$C_{q}\not \in {\cal M}$. Let  $p\in \pi$. Then $C_{q}$  is 
$M(p)$-critical and so  $C_{q} \in {\cal F}\subseteq {\cal M}$ by the hypothesis. This 
contradiction shows that   ${\cal N}\subseteq {\cal 
M}$.

 Now we show that  $\text{Int}_{\cal  F}(G)\leq
Z_{\pi{\cal M}}(G)$ for every group $G$.
Suppose that this assertion   is false  and let $G$ be a
counterexample of minimal order. Let $I= \text{Int}_{\cal F}(G)$
and   $Z=Z_{\pi \cal M}(G)$.  Then  $1 < I < G$ and  $Z \ne G$. Let $N$
be a minimal normal subgroup of $G$ and  $L$  a minimal normal
subgroup of $G$ contained  in $ I$. Then $\pi (L)\leq \pi ({\cal
F})$. We proceed via the following steps.

(1)    $IN/N\leq  \text{Int}_{\cal F}(G/N)\leq Z_{\pi \cal
M}(G/N)$.

Indeed, by Lemma 2.4(a), we have $IN/N\leq  \text{Int}_{\cal
F}(G/N)$. On the other hand, by the choice of $G$, $\text{Int}_{\cal
F}(G/N)\leq Z_{\pi{\cal M}}(G/N)$.

(2)   {\sl $L\nleq Z$; in particular, the order of $L$ is divisible
by some prime $p\in \pi$}.

Suppose that $L\leq Z$. Then, clearly,  $Z/L=Z_{\pi{\cal M}}(G/L)$,
and $I/L= \text{Int}_{\cal F}(G/L)$ by  Lemma 2.4(e).  But by (1),
$\text{Int}_{\cal F}(G/L) \leq Z_{\pi{\cal M}}(G/L)$. Hence $I/L\leq
Z/L$. Consequently, $I\leq Z$, a contradiction.

(3) {\sl If $L\leq M < G$, then $L\leq Z_{{\cal M}}(M)$}.

By Lemma 2.4(b), $L\leq I\cap M\leq \text{Int}_{\cal F}(M)$. But
since $|M| < |G|$, we have that $ \text{Int}_{\cal F}(M)\leq
Z_{\pi{\cal M}}(M)$ by the choice of $G$.  Hence  $L\leq Z_{\pi{\cal
M}}(M)$ and so $L\leq Z_{\cal M}(M)$ since the order of $L$ is
divisible by some prime $p\in \pi$ by (2).

(4) {\sl $L=N$ is the  unique minimal normal subgroup of $G$}.

Suppose that $L\ne N$. Then by (1),  $NL/N\leq Z_{\pi{\cal
M}}(G/N)$. Hence from the $G$-isomorphism $NL/N\simeq L$ we obtain
$L\leq Z$, which contradicts (2).

(5) $L\nleq \Phi (G)$.

Suppose that $L\leq \Phi (G)$. Then $L$ is a $p$-group by (2). Let
$C=C_{G}(L)$ and $M$ be any maximal subgroup of $G$. Then $L\leq M$.
Hence  $L\leq Z_{\cal M}(M)$ by (3), so $M/M\cap C\in M(p)$ by
Lemmas 2.1(1) and Lemma 2.5. If  $C\nleq M$, then  $G/C=CM/C\simeq
M/M\cap C\in M(p)$. This implies that $L\leq Z$, which contradicts
(2). Hence $C\leq M$ for all maximal subgroups $M$ of $G$. It
follows that $C$ is nilpotent. Then in view of (4), $C$ is a
$p$-group since $C$ is normal in $G$. Hence for every maximal
subgroup $M$ of $G$
 we have $M\in {\cal G}_{p}M(p)=M(p)$. But  since  $M(p)\subseteq {\cal
M}$,  $G\not \in M(p)$ (otherwise $G\in {\cal M}$ and so $G=Z$).
This shows that $G$ is an $M(p)$-critical group. Therefore $G\in
{\cal F}$ by the hypothesis. But since ${\cal F}\subseteq {\cal M}$,
we have $G\in {\cal M}$ and so $G=Z$, a contradiction. Thus (5)
holds.

(6) {\sl $C=C_{G}(L)\leq L$} (This follows from (4), (5) and Theorem
15.6 in \cite[Chapter A]{DH}).

(7) {\sl If  $L\leq M < G$, then  $M\in M(p)$.}

First by (3), $L\leq Z_{{\cal M}}(M)$. If $L=C$, then $M/L=M/M\cap
C\in M(p)$ by Lemma 2.5, which implies that $M\in {\cal
G}_{p}M(p)=M(p)$ since $L$ is a $p$-group by (2). Now suppose that
$L$ is a non-abelian group. Let $1=L_{0} < L_{1}  < \ldots < L_{n}=L
$ be a chief series of $M$ below $L$. Let
$C_{i}=C_{M}(L_{i}/L_{i-1})$ and $C_{0}=C_{1}\cap \ldots \cap
C_{n}$. Since $L\leq Z_{\cal M}(M),$  $M/C_{i}\in M(p)$ for all
$i=1, \ldots , n $. It follows that $M/C_{0}\in M(p)$. Since $C=1$
by (4) and (6), for any minimal normal subgroup $R$ of $M$ we have
$R\leq L$. Suppose that $C_0\ne 1$ and let $R$ be a minimal normal
subgroup of $M$ contained in $C_0$. Then  $R\leq L$ and  $R\leq
C_M(H/K)$ for each  chief factor $H/K$  of $M$. Thus $R\leq F(M)$ is
abelian and hence $L$ is abelian. This contradiction shows that
 $C_{0}=1$. Consequently, $M\in M(p)$.

(8) {\sl  There exists a subgroup $U$ of $G$ such that $U\in {\cal
F}$ and $LU=G$.}

Indeed, suppose that  every maximal subgroup of $G$ not containing $L$ belongs
to $M(p)$. Then by (7), $G$ is an $M(p)$-critical group. Hence $G\in
{\cal F}$ by the hypothesis. But then $I=G$, a contradiction.  Hence there
 exists a maximal subgroup $M$ of $G$ such
that $G=LM$ and $M\not\in M(p).$ Take an $M(p)$-critical subgroup
$U$ of $M$. Then in view of (7), $LU=G$ and $U\in {\cal F}$ by the
hypothesis.

(9) {\sl The final contradiction for (b).}

Since $L\leq I$ and $G/L=UL/L\simeq U/U\cap L\in {\cal F}$ by (8),
it follows from Lemma 2.4(c) that $G\in {\cal F}$ and so $G=I$. The
final contradiction shows that $\text{Int}_{\cal  F}(G)\leq
Z_{\pi{\cal M}}(G)$ for every group $G$.   The second assertion of (b) can 
be proved similarly.  The theorem is proved.

{\bf Proofs of Theorems A and B.}  Since $Z_{\cal F}(G)\leq
Int_{\cal F}(G)$ by Lemma 2.4(g), the sufficiency is a special case,
when ${\cal F}={\cal M}$, of Theorem C (b). Now suppose that the
equality $Z_{\pi{\cal F}}(G)= \text{Int}_{\cal F}(G)$ holds for each
 (soluble) group $G$.

First we show that  ${\cal N}\subseteq {\cal F}$. Let  $F$ be the  canonical local
 satellite   of $ {\cal F}$. Suppose that for some 
group $C_{q}$  of prime order $q$ we have $C_{q}\not \in {\cal F}$. Let 
$p\in \pi $ and $G=PC_{q}$, where $P$ is  a
simple ${\mathbb F}_{p}A$-module $P$    which is     faithful  for
$C_{q}$. Then $P=\text{Int}_{\cal F}(G)$ and $Z_{\cal F}(G)=1$ since 
$F(p)\subseteq {\cal F}$. This contradiction shows that  ${\cal 
N}\subseteq {\cal F}$.

Now  we  show that ${\cal G}_{\pi'}{\cal F}={\cal F}$ (${\cal 
S}_{\pi'}{\cal F}={\cal F}$, respectively). The
inclusion ${\cal F}\subseteq {\cal S}_{\pi'}{\cal F}$ is evident.
Suppose that ${\cal G}_{\pi'}{\cal F}\not\subseteq  {\cal F}$ 
  (${\cal S}_{\pi'}{\cal F}\not\subseteq  {\cal F}$) and
let $G$ be a group of minimal order in ${\cal G}_{\pi'}{\cal F} \setminus  {\cal F}$
 (in ${\cal S}_{\pi'}{\cal F} \setminus  {\cal F}$, respectively).
 Then  $ G^{\cal F}$  is the unique minimal
normal subgroup of $G$ and  $G^{\cal F}$ is a $\pi'$-group. Hence
$$G^{\cal F}\leq Z_{\pi{\cal F}}(G)= \text{Int}_{\cal F}(G).$$
It follows from Lemma 2.4(c) that $G\in {\cal F}$. This
contradiction shows that ${\cal G}_{\pi'}{\cal F}={\cal F}$ 
(${\cal S}_{\pi'}{\cal F}={\cal F}$, respectively).

Finally, we  show    that   $\cal F$  satisfies the $\pi$-boundary
condition (the $\pi$-boundary
condition in the class $\cal S$, respectively). Suppose that this is
false.  Then  for some $p\in \pi$, the set of all (soluble)
$F(p)$-critical groups $A$ with $A\not \in {\cal F}$  is non-empty.
Let us choose in this set a group $A$ with
 minimal $|A|$.
 Then  by Lemma 2.6,
$A^{\cal F}$ is the unique minimal normal subgroup of $G$ and
$O_{p}(A)=1=\Phi (A)$. Hence by  Lemma 2.10, there
exists  a  simple ${\mathbb F}_{p}A$-module  $P$    which is
faithful  for $A$. Let $G=P\rtimes A$ and  $M$ be any maximal
subgroup of $G$. If $P\nleq M$, then $M\simeq G/P\simeq A \not \in
{\cal F}$. On the other hand, if $P\leq M$, then $M=M\cap PA=P(M\cap
A)$, where $M\cap A$ is a maximal subgroup of $A$. Hence   $M\cap
A\in F(p)$ and so  $M\in {\cal G}_{p}F(p)=F(p)\subseteq {\cal F}$.
Therefore $P$ is contained in the intersection of all $\cal
F$-maximal subgroups of $G$. Then $P\leq Z_{\pi{\cal F}}(G)$ by our
assumption on $\cal F$. It follows that $A\simeq G/P=G/C_{G}(P)\in
F(p)\subseteq {\cal F}$ by Lemmas 2.2(1) and Lemma 2.5. But this
contradicts the choice of $A$. Therefore  $\cal F$
 satisfies the $\pi$-boundary condition ($\cal F$  satisfies the $\pi$-boundary
condition  in the class $\cal S$).  The theorems are proved.

In view of Theorems A, B   and C we have 

{\bf Corollary 3.1}.  {\sl   Let
 ${\cal F}$  be a  hereditary saturated formation
 with   $(1) \ne {\cal F}\ne {\cal G}$. Then the equality
  $ \text{Int}_{\cal F}(G)=Z_{\cal F}(G) $  holds for each
 group $G$ if and only if ${\cal N}\subseteq {\cal F}$  and $\cal F$  satisfies
 the boundary condition. }

{ \bf Corollary 3.2}. {\sl    Let
${\cal F}$  be a  hereditary saturated formation of soluble groups
 with  $(1) \ne {\cal F}\ne {\cal S}$. Then the  equality  $ \text{Int}_{\cal F}(G)=Z_{\cal
F}(G) $  holds for each soluble group $G$ if and only if 
 ${\cal N}\subseteq {\cal F}$,  ${\cal F}$     
satisfies the boundary condition    in the class ${\cal S}$. }

Note that Corollary 3.2
 also follows directly from \cite[Main Theorem]{BeidH}.

{\bf Proof of  Theorem  D.} We will prove the theorem by induction
on $|G|$. If $G \in {\cal F}$, then  $$\text{Int}^{*}_{\cal
F}(G)=G=\text{Int}_{\cal F}(G).$$ We may, therefore, assume that $G
\not\in {\cal F}$. Let $I=\text{Int}_{\cal F}(G)$,
$I^{*}=\text{Int}^{*}_{\cal F}(G)$ and $N$ a minimal normal subgroup
of $G$. Then $I \leq I^{*}$. Hence we may assume that $I^{*}\ne 1.$

(1) {\sl $I^{*}N/N\leq \text{Int}^{*}_{\cal F}(G/N)$. }

If $U/  N$ is a non-$K$-${\cal F}$-subnormal  $\cal F$-maximal
subgroup of $G/N$, then for some non-$K$-${\cal F}$-subnormal $\cal
F$-maximal subgroup $U_{0}$ of $G$ we have  $U=U_{0}N$ by Lemma
2.8(i). Let $$\text{Int}^{*}_{\cal F}(G/N)= U_{1}/N
  \cap  \ldots
\cap U_{t}/N, $$ where  $U_{i}/N$ is a non-$K$-${\cal F}$-subnormal
$\cal F$-maximal subgroup of $G/N$ for all $i=1,  \ldots , t$. Let
$V_{i}$ be a non-$K$-${\cal F}$-subnormal  $\cal F$-maximal subgroup
of $G$ such that  $U_{i}=V_{i}N$. Then $I^{*}\leq  V_{1}  \cap
\ldots \cap V_{t}$. Hence $I^*N/N\leq \text{Int}^{*}_{\cal F}(G/N)$.

(2) {\sl If $N\leq I^{*}$, then $\text{Int}^{*}_{\cal F}(G/N)=I^{*}/N$.}

By Lemma 2.8(i), it is enough to prove that if  $U$ is a
non-$K$-${\cal F}$-subnormal $\cal F$-maximal subgroup of $G$, then
$U/N$  is a non-$K$-${\cal F}$-subnormal $\cal F$-maximal subgroup
of $G/N$. Let $U/N\leq X/N$, where $X/N$  is a  non-$K$-${\cal
F}$-subnormal
 $\cal F$-maximal subgroup of $G/N$. By Lemma 2.8(i),  $X=U_{0}N$
for some  non-$K$-${\cal F}$-subnormal $\cal F$-maximal subgroup
$U_{0}$ of $G$. But since $N\leq U_{0}$,   $U/N\leq U_{0}/N$ and so
$U=U_{0}$. Thus $U/N= X/N$.

(3) {\sl If $I^{*}\cap H\leq \text{Int}^{*}_{\cal F}(H)$ for any
subgroup $H$ of $G$.}

Let $V$ be an arbitrary  non-$K$-${\cal F}$-subnormal  $\cal
F$-maximal subgroup of $H$. Then $V=H\cap U$ for some non-$K$-${\cal
F}$-subnormal $\cal F$-maximal   subgroup $U$ of $G$ by Lemma
2.8(ii). Thus there are non-$K$-${\cal F}$-subnormal $\cal
F$-maximal subgroups $U_{1},   \ldots , U_{t}$ of $G$ such that
$$\text{Int}^{*}_{\cal F}(H)=U_{1} \cap   \ldots \cap U_{t}\cap H.$$
This induces that  $I^{*}\cap H\leq  \text{Int}^{*}_{\cal F}(H)$.

(4) {\sl If $E=RV$ for some    normal subgroup $R$ of $G$ contained in $I^{*}$ and 
 $K$-${\cal F}$-subnormal   subgroup  $V\in {\cal
F}$, then $E\in {\cal F}$.}

First note that    by (3), $R\leq
\text{Int}^{*}_{\cal F}(E)$. On the other hand, by Lemma 2.7(i), $V$
is a $K$-${\cal F}$-subnormal   subgroup of $E$.   Hence we
 need only consider the case when $G=E$. Assume that 
 $G\not\in {\cal F}$.   Then 
 $D\ne 1$.  Let $R$ be any minimal normal
subgroup of $G$. Then $(DR/R)(VR/R)=G/R$, where $DR/R\leq
\text{Int}^{*}_{\cal F}(G/R)$  by (2),  and $VR/R\simeq V/V\cap R\in {\cal F}$.
On the other hand, by induction we have $\text{Int}^{*}_{\cal 
F}(G/R)=\text{Int}_{\cal 
F}(G/R)$, so $G/R\in {\cal F}$ by Lemma 2.4 (d).    This
implies that $R$ is the only minimal normal subgroup of $G$ and so
$R=G^{\cal F}\leq \text{Int}^{*}_{\cal F}(G)$.
 Let $W$
be a minimal supplement of $R$ in $G$. Then $W\in {\cal F}$ by Lemma
2.3(i). Let $W\leq E$, where $E$ is an ${\cal F}$-maximal subgroup
of $G$. If $E$ is not $K$-${\cal F}$-subnormal in $G$, then $R
\leq
E$. Thus $G=RW=RE=E\in {\cal F}$, a contradiction. This shows that
$E$ is $K$-${\cal F}$-subnormal in $G$. But then there is a proper
subgroup $X$ of $G$ such that $E\leq X$ and either $X$ is normal in
$G$ or $R=G^{\cal F}\leq X$. In  both of this cases, we have that
$G=RV=RX=X < G$, a contradiction. Hence we have  (4).

{\sl Conclusion.}

Let $R$ be a minimal normal subgroup
of $G$ contained in $I^{*}$. If $R\leq I$, then
$I/R=\text{Int}_{\cal F}(G/R)$ by Lemma 2.3(e), and
$I^{*}/R=\text{Int}^{*}_{\cal F}(G/R)$ by (2). Therefore by
induction, $\text{Int}^{*}_{\cal F}(G/R)=\text{Int}_{\cal F}(G/R)$.
It follows that $I=I^{*}$.

Finally, suppose that $R\nleq I$. Then $R\nleq U$ for some $\cal
F$-maximal  subgroup $U$ of  $G$.  It is clear that $U$
is a $K$-${\cal F}$-subnormal   subgroup of  $G$ and hence  $E\in {\cal
F}$ by (4). But then  $E=U$, which implies $R\leq U$, a
contradiction. The theorem is proved.

\section{Applications and Remarks}

{\bf Applications of Theorems A, B and D.} We say  that  $\cal F$
\emph{satisfies the $p$-boundary condition}
 if $\cal F$ satisfies the  
$\{p\}$-boundary condition
in the class of all groups.

{\bf Lemma 4.1.} {\sl
Let  ${\cal F}=LF(F)$, where  $F$  the canonical local satellite of ${\cal F}$.
  Suppose that for some prime $p$ we have $F(p)
={\cal F}$. Then ${\cal F}$  does not satisfy the $p$-boundary condition.}

{\bf Proof.}  Indeed, in this case  every     ${\cal F}$-critical
group is also  $F(p)$-critical.

A group $G$ is called $\pi$-closed if $G$ has a normal Hall $\pi$-subgroup.

{\bf Proposition 4.2.} {\sl The formation ${\cal F}$  of all
$\pi$-closed groups satisfies the $\pi'$-boundary condition, but
${\cal F}$  does not satisfy the $p$-boundary condition for any
$p\in \pi$. }

{\bf Proof. }   Let ${\cal F}={\cal G}_{\pi}{\cal G}_{\pi'}$ be  the
formation of all $\pi$-closed  groups, $F$ the canonical local
satellite of ${\cal F}$. Then $F(p)={\cal F}$  for all $p\in \pi$,
and  $F(p)={\cal G}_{\pi'}$ for all primes $p\in \pi'$ by Theorem
3.1.20 in \cite {Guo}. Hence ${\cal F}$ satisfies the
$\pi'$-boundary condition and does not satisfy the $p$-boundary
condition for any $p\in \pi$ by Lemma  4.1.

{\bf Lemma  4.3.} {\sl Let $\{{\cal F}_{i} \mid i\in I \}$ be any set
of non-empty saturated formations and ${\cal F}= \cap _{i\in I}
{\cal F}_{i}$.}

(1) {\sl If for each $i\in I$, ${\cal F}_{i}$   satisfies the
$p$-boundary condition, then $\cal F$ satisfies the $p$-boundary
condition.}

(2) {\sl Suppose that $I= \{1, 2 \}$,  $F_{i}$ is the canonical
local satellite of ${\cal F}_{i}$ and that there is a set $\pi$ of
primes satisfying the following conditions:}

(a) {\sl ${\cal F}_{1}$ satisfies the  $\pi$-boundary condition, and
for any $p\in \pi$, we have $F_{1}(p)\subseteq {\cal F}_{2}=
F_{2}(p)$  and every $F_{1}(p)$-critical group belongs to ${\cal
F}_{2}$ }.

(b) {\sl ${\cal F}_{2}$ satisfies the $\pi'$-boundary condition, and
for any $p\in \pi'$, we have $F_{2}(p)\subseteq {\cal F}_{1}=
F_{1}(p)$  and every $F_{2}(p)$-critical group belongs to ${\cal
F}_{1}$ }.

{\sl Then ${\cal F}$ satisfies the boundary condition.}

{\bf Proof.}  (1) Let $F_{i}$ be the canonical local satellite of
${\cal F}_{i}$ and $F$  the canonical local satellite of ${\cal F}$.
If $f(p)=\cap _{i\in I} F_{i}(p)$, then  $F(p)={\cal G}_{p}f(p)$ by
Theorem 3.3 in \cite[Chapter 1]{26}. Now let $G$ be any
$F(p)$-critical group, $i\in I$. Since $F(p)\subseteq F_{i}(p)$, all
maximal subgroup of $G$ belongs to $F_{i}(p)$. Hence $G\in {\cal
F}_{i}$ since $F_{i}(p)\subseteq {\cal F}_{i}$ and ${\cal F}_{i}$
satisfies the $p$-boundary condition. This implies that $G\in {\cal
F}$ and therefore $\cal F$ satisfies the $p$-boundary condition.

(2) In this case,  $F(p)= F_{1}(p)$ for all $p\in \pi$ and  $F(p)=
F_{2}(p)$ for all $p\in \pi'$, where $F$ is the canonical local
satellite of ${\cal F}$. Hence if $p\in \pi$ and $G$ is an
$F(p)$-critical groups, then $G\in {\cal F}$ by hypothesis (a). This
shows that ${\cal F}$ satisfies the $\pi$-boundary condition.
Similarly we see that ${\cal F}$ satisfies the $\pi'$-boundary
condition.

A group $G$ is called  $p$-decomposable if there exists a subgroup $H$ of
$G$ such that $G=P\times H$ for some (and hence the unique)
 Sylow $p$-subgroup $P$ of $G$.

{\bf Corollary  4.4.} {\sl The formation of all $p$-decomposable
groups satisfies the boundary condition. }

{\bf Proof. } Let $\cal F$ be the formation of all $p$-decomposable
groups. Then ${\cal F}={\cal G}_{p'}{\cal G}_{p}\cap {\cal
G}_{p}{\cal G}_{p'}$. Hence the assertion follows from Proposition
4.2 and   Lemma 4.3(2).

From Corollary 4.4 and Theorem A we get

{\bf Corollary 4.5.} {\sl Let $D$ be the intersection of all maximal
$p$-decomposable  subgroups of $G$. Then $D$  is the largest normal
subgroup of $G$ satisfying $D=O_{p'}(D)\times O_{p}(D)$, and $G$
induces the trivial automorphisms group on every chief factor of $G$
below $O_{p}(D)$ and a   $\pi'$-group of automorphisms   on every
chief factor of $G$ below $O_{p'}(D)$.}

Since a $p$-nilpotent group is $p'$-closed, the following result
directly follows from Proposition 4.2.

{\bf Corollary  4.6.}  {\sl The formation of all $p$-nilpotent
groups satisfies the $p$-boundary condition. }

From Corollary 4.6 and Theorem A, we have

{\bf Corollary 4.7.}  {\sl Let $D$ be the intersection of all
maximal $p$-nilpotent subgroups of a group $G$. Then $D$  is the
largest normal subgroup of $G$ satisfying $O_{p'}(D)= O_{p'}(G)$,
and $D/O_{p'}(G)\leq Z_{\infty}(G/O_{p'}(G))$.}

Note that another proof of Corollary 4.7 was obtained in the paper
\cite{BeidH}.

{\bf Remark 4.8.} If $\pi \ne \{ 2 \}$, then the formation $\cal F$
of all $\pi$-supersoluble groups does not satisfy the $\pi$-boundary
condition. Indeed, let $F$ be the canonical local satellite of
${\cal F}$. Then $F(p)={\cal N}_{p}{\cal A}(p-1)$,  where ${\cal
A}(p-1)$ is the formation of all abelian groups of exponent dividing
$p-1$ \cite[p.358]{DH}, for all $p\in \pi $ and $F(p)={\cal F}$ for
all primes $q\not \in \pi  $ (see Example 3.4(e) and Theorem 4.8 in
\cite[Chapter IV]{DH}). Let $2\ne p \in \pi $ and   $q$  any prime
with $q$ divides $p-1$. Let $G=Q\rtimes C$, where   $C$ is a group
of order $p$ and $Q$ is an $\Bbb{F}_{q}C$-module which is faithful
for $C$. Then $G$ is a non-supersoluble $F(p)$-critical groups.

{\bf Proposition 4.9.}  {\sl Let ${\cal F}={\cal N}{\cal L}$.}

(i) {\sl If ${\cal L}$ is a hereditary saturated formation   satisfying 
the boundary condition    in the class  of all
soluble   groups, then ${\cal F}$ satisfies the
boundary condition  in the class  of all
soluble   groups. }

(ii) {\sl If ${\cal L}$ is a formation of nilpotent  groups
with $\pi ({\cal L})=\Bbb{P}$, then ${\cal F}$ satisfies the
boundary condition. }

{\bf Proof.} Let $F$ be the canonical local satellite of ${\cal F}$.
Then by Lemma 2.9(2),  $F(p)={\cal G}_{p}{\cal L}$ for
 all primes $p$. Assume that  ${\cal F}$ does not   satisfy
the boundary condition (${\cal F}$ does not   satisfy
the boundary condition in the class  of all
soluble   groups). Then for some
prime $p$, the set of all $F(p)$-critical (soluble) groups  $A$ with  $A\not \in {\cal F}$
 is non-empty. Let $G $ be a group of minimal order in this set.  Then 
 $L=G^{\cal F}$ is the unique minimal normal
 subgroup of $G$ and 
  $O_{p}(G)=1=\Phi (G)$  by Lemma 2.6.   First 
suppose that $G$ is soluble. Then   $L=C_{G}(L)$ is a $q$-group for some prime $q\ne p$
and $G=L\rtimes M$ for some maximal subgroup $M$  of $G$ with $O_{q}(M)=1$ 
by  
\cite[Chapter A, Theorem 15.6]{DH}. Let $M_{1}$ be any  maximal
subgroup of $M$. Then  $LM_{1}\in F(p)$, so $LM_{1}\in {\cal H}$  since $L=C_{G}(L)$. 

(i) Let $L$ be the canonical local satellite of ${\cal L}$. Since $LM_{1}\in {\cal L}$,
 $$LM_{1}/O_{q', q}(LM_{1})=LM_{1}/O_{q}(LM_{1})=LM_{1}/LO_{q}(M_{1})\simeq M_{1}/M_{1}\cap
 LO_{q}(M_{1} \in L(q)={\cal G}_{p}L(q).$$   Hence every  maximal subgroup 
of $M$ belongs to $L(q)$. Since ${\cal L}$   satisfies 
the boundary condition    in the class  of all
soluble   groups, it follows that $G\in   {\cal F}$. This contradiction 
completes the poof of (i).

(ii) We may suppose that $M_{1}$ is  a normal maximal
subgroup of $M$. Then  $LM_{1}\in F(p)$, so $LM_{1}\in {\cal H}$ 
since $L=C_{G}(L)$. This implies that $M_{1}=1$. Consequently
$|M|$ is prime, so $M\in  {\cal L}$ since $\pi ({\cal L})=\Bbb{P}$. But 
then $G\in {\cal F}={\cal N}{\cal L}$.  Therefore $G$ is not soluble.
 Let
$q\ne p$ be any prime divisor of $|G|$. Suppose that $G$ is not
$q$-nilpotent. Then $G$ has a $q$-closed $\cal N$-critical  subgroup
$H=Q\rtimes R$ by \cite[IV, Theorem 5.4]{Hupp}, where $Q$ is a Sylow
$q$-subgroup of $H$, $R$ is a cyclic Sylow $r$-subgroup of $H$.
Since $G$ is not soluble, $H\ne G$. Hence  $H\leq M\in F(p)$ for
some  maximal subgroup $M$ of $G$. Since $M\in {\cal G}_{p}{\cal
L
}$,  $M^{\cal N}\leq O_{p}(M)$ and hence  $H^{\cal N}\leq Q\cap O_{p}(H)=1$.
Therefore $H$ is nilpotent. This contradiction shows that  $G$ is
$q$-nilpotent for all primes $q\ne p$. This induces that $G^{\cal
N}$ is a Sylow $p$-subgroup of $G$ and thereby $G$ is soluble. This contradiction 
completes the poof of (ii).

{\bf Remark 4.10.}  The condition "${\cal L}$ is a hereditary saturated formation   satisfying 
the boundary condition    in the class  of all
soluble   groups" can not be omitted in the Statement (ii). Indeed, let
 ${\cal F}={\cal N}{\cal U}$ and  $G=P\rtimes A_{4}$, where $P$ is
 a simple ${\mathbb F}_{3}A_{4}$-module $P$    which is     faithful  for
$A_{4}$. Let $F$ be the canonical local satellite of  ${\cal F}$. Then 
$F(2)={\cal F}_{2}{\cal U}$ by Lemma 2.9 (2).  Therefore $G$ is an  
$F(2)$-critical group and   $G\not \in {\cal F}$. Thus $G$ does not 
satisfy  the boundary condition    in the class  of all
soluble   groups.

{\bf Corollary 4.11.} {\sl   Let ${\cal F}$  be the class of all
 groups with  $G'\leq F(G)$. Then ${\cal F}$  satisfies the boundary condition. }

From Corollary 4.11 and Theorem A, we obtain

{\bf Corollary 4.12.} {\sl Let $D$ be the intersection of all
maximal subgroups $H$ of $G$ with the property $H'\leq F(H)$. Then $D$
is the largest normal subgroup of $G$ such that $D'\leq F(D)$ and
$G$ induces  an abelian group of  automorphisms on every chief
factor of $G$ below $D$.}

Note that Corollary 4.12 can be also found in \cite[Corollary 7 and Remark 4]{BeidH}.

Following \cite[Chapter VII, Definitions 6.9]{DH} we write  $l(G)$
to denote the nilpotent length of the group $G$.  Recall that ${\cal
N}^{r}$ is the product of $r$ copies of ${\cal N}$; ${\cal N}^{0}$
is the class of groups of order 1 by definition. It is well known
that ${\cal N}^{r}$ is the class of all soluble groups $G$ with
$l(G) \leq r$. It is also known that  ${\cal N}^{r}$ is a hereditary
saturated formation (see, for example,  \cite[p. 358]{DH}). 
Hence from Proposition 4.9  we get.

{\bf Corollary  4.13.} {\sl Let ${\cal F}={\cal N}^{r}{\cal L}$ ($r\in
{\mathbb N})$,
where ${\cal L}$ is a subformation of the formation of all abelian
groups with $\pi (L)=\Bbb{P}.$  Then ${\cal F}$ satisfies the
boundary condition in the class of all soluble groups. }

From Proposition  4.9 and Theorem A we get

{\bf Corollary 4.14.} {\sl Let $D$ be the intersection of all
maximal metanilpotent  subgroups of a  group $G$. Then $D$ is the largest
normal subgroup of $G$ satisfying $D$ is metanilpotent  and $G$
induces a nilpotent  group of  automorphisms on every chief factor
of $G$ below $D$.}

It is clear that every subnormal subgroup is a $K$-${\cal
F}$-subnormal    subgroup as well. On the other hand, in the case, where ${\cal N}\subseteq
 {\cal F}$,  every  $K$-${\cal
F}$-subnormal    subgroup of any soluble subgroup $G$ is  ${\cal
F}$-subnormal in $G$.    Hence from Theorems D and the
above corollaries we get

{\bf Corollary 4.15.} {\sl Let  ${\cal F}$  be  the class of all groups $G$
 with $G'\leq F(G)$. }

(i) {\sl The subgroup $ Z_{\cal F}(G)$ may be characterized as the
intersection of all non-subnormal  $\cal F$-maximal  subgroups of
$G$, for each  group $G$. }

(ii) {\sl The subgroup $ Z_{\cal F}(G)$ may be characterized as the
intersection of all non- ${\cal
F}$-subnormal   $\cal F$-maximal  subgroups of
$G$, for each  soluble group $G$. }

{\bf Corollary 4.16.} {\sl Let  ${\cal F}$  be one of the following
formations: }

(1) {\sl the class of all nilpotent groups} (Baer \cite{BaerI});

(2) {\sl the class of all groups $G$ with $G'\leq F(G)$} (Skiba \cite{skjpaa});

(3) {\sl the class of all $p$-decomposable groups}  (Skiba
\cite{Skibaja5}).

{\sl Then for each  group $G$,   the equality }
 $ \text{Int}_{\cal F}(G)=Z_{\cal F}(G)$ {\sl holds.}

{\bf Corollary 4.17} (Sidorov \cite{Sidorov}).  {\sl Let  ${\cal F}
$ be  the class of all soluble  groups $G$ with $l(G) \leq r \ (r\in
{\mathbb N})$}. Then for each  soluble group $G$, the equality
$Z_{\cal F}(G)=\text{Int}_{\cal F}(G)$ {\sl holds.}

{\bf An application of Theorem C (a).}  Let $p_{1} > p_{2} > \ldots  > p_{r}$  
   be the distinct primes dividing $|G|$, $P_{i}$  a Sylow 
$p_{i}$-subgroup of $G$. Then $G$ is said to satisfy \emph{the Sylow tower  
property} if all subgroups  $P_{1}$, $P_{1}P_{2}$,  \ldots   , $P_{1}P_{2}\ldots P_{r-1}$  
 are normal in $G$.  Every supersoluble group  satisfy the Sylow tower  
property.

{\bf Proposition 4.18.} Let ${\cal F}$ be the formation of all
supersoluble groups and ${\cal M}$  be the formation of all
  groups satisfying the Sylow tower  property. Then $\text{Int}_{\cal F}(G)  \leq
\text{Int}_{\cal M}(G)$ for any  group $G$.

{\bf Proof.} Let $F$ be the canonical local satellite of ${\cal F}$.
Then $F(p)={\cal G}_{p}{\cal A}(p-1)$ for  primes $p$ (see Remark 4.8).
 Let $G$ be any  $F(p)$-critical group satisfying the Sylow tower  
property.
We show that $G\in {\cal F}$. Let $q$ be the largest prime dividing
$|G|$, $P$ the Sylow $q$-subgroup of $G$. If $G=P$, then clearly
$G\in {\cal F}$. Let $P\ne G$. Then every Sylow subgroup of $G$
belongs to $F(p)$. Hence $q=p$ and if $E$ is a Hall $p'$-subgroup of
$G$, then $E\in {\cal A}(p-1)$. But then $G\in F(p)$, a
contradiction. Therefore $G\in \cal F$. This shows that ${\cal F}$
satisfies the boundary condition in ${\cal M}$. Hence, in view of
Theorem C(b) , we have $\text{Int}_{\cal F}(G) \leq \text{Int}_{\cal
M}(G)$.

{\bf Remark 4.19.} If ${\cal F}\subseteq {\cal M}$ be hereditary
saturated formations, then for every group $G$ we have $Z_{\cal
F}(G)\leq Z_{\cal M}(G)$ and $Z_{\cal F}(G) \leq  \text{Int}_{\cal
F}(G)$ by Lemma 2.4(g). Therefore, if ${\cal F}$
satisfies the boundary condition, then $\text{Int}_{\cal F}(G) \leq
\text{Int}_{\cal M}(G)$ for every group $G$. But, by Remark 4.8, 
 ${\cal U}$ does not necessarily satisfy the boundary
condition.  Hence we can not deduce  Proposition 4.18 from  Theorem
A.

{\bf Lemma 4.20.} {\sl Let  ${\cal F}=LF(F)$
be a  saturated formations, where $F$  is the canonical local satellites of ${\cal
F}$. Let ${\cal M}= {\cal F}\cap {\cal S}$ and $M$   the canonical local satellites of ${\cal
M}$. Then $M(p)= F(p)\cap {\cal S}$ for all primes $p$.}

{\bf Proof.}  It is clear that  ${\cal M}=LF(F_{1})$, where  $F_{1}(p)= 
F(p)\cap {\cal S}$  for all primes $p$. On the other  hand, for any $p$ we 
have  $F_{1}(p)\subseteq  {\cal M}$,  and  $F_{1}(p)={\cal G}_{p}F_{1}(p)$ 
since  $${\cal G}_{p}F_{1}(p)={\cal G}_{p}(F(p)\cap {\cal S})\subseteq 
F(p)\cap {\cal S}.$$   Hence $F_{1}=M$ is the canonical local satellite of 
${\cal M}$.

{\bf Lemma 4.21.} {\sl Let  $(1)\ne  {\cal F}=LF(F)$
be a  saturated formations  and  ${\cal M}= {\cal F}\cap {\cal S}$. Let $\pi \subseteq
 \pi ({\cal F})$. Then ${\cal F}$ satisfies the boundary condition    in the class $\cal S$
 of all
soluble   groups if and only if ${\cal M}$ satisfies the boundary condition  
  in  $\cal S$.}

 {\bf Proof.} Let  $F$  and $ M$ be the canonical  local  satellite of the formation
 ${\cal F}$  and  ${\cal M}$,  respectively. Let $p\in \pi$. Then  $M(p) = F(p)\cap {\cal S}$
by Lemma 4.20. Hence a soluble group 
$G$ is   $F(p)$-critical group   if and only if  $G$ is   $M(p)$-critical. 
On the other hand, the group $G$ belongs to ${\cal F}$ if and only if  $G$ 
belongs to ${\cal M}$.

In view of Lemmas 4.20 and 4.21, 
  Theorem B may be proved in the following general  form.

{\bf Theorem 4.22.}   {\sl  Let ${\cal F}$  be a   saturated
formation such that  ${\cal M}={\cal F}\cap {\cal S}$ is hereditary and 
   $(1) \ne {\cal M} \ne {\cal S}$.  Let $ \pi\subseteq
\pi (\cal F)$.  Then  the equality
$$Z_{\pi \cal F}(G)= \text{Int}_{\cal F}(G)$$  holds for each soluble
group $G$ if and only if ${\cal N}\subseteq {\cal F}={\cal S}_{\pi'}  
  {\cal F}$ and
$\cal F$  satisfies the boundary condition in the class  of all
soluble   groups. }

\end{document}